\documentclass{iise}% use this for final manuscript

\usepackage{multirow}

\conference{Proceedings of the IISE Annual Conference \& Expo 2026\\
Y. Xiang, D., Yu, R. Thiesing, eds.}% Do not change this line.

\title{\titlesize Distributed and Dynamic Hub Network Operation Planning in a Hyperconnected Less-Than-Truckload Operating System}

\author{
Tiankuo Zhang, Jihye Jung, Paria Nourmohammadi, Benoit Montreuil, Alan Erera, Sahrish Jaleel Shaikh\\Physical Internet Center, Supply Chain \& Logistics Institute, School of Industrial \& Systems Engineering, Georgia Institute of Technology
\\Atlanta, GA}
\authorlist{Zhang, Jung, Nourmohammadi, Montreuil, Erera, and Shaikh}%Heading
\abstractID{16686}% Fill in the web conference management system-assigned abstract ID.

\begin{document}
\maketitle

\begin{abstract}
{\small The less-than-truckload (LTL) industry plays a vital role in enhancing the efficiency and sustainability of logistics systems, as LTL shipments offer greater consolidation opportunities than full-truckload shipments. Despite of this flexibility, the average cost of LTL shipments remains considerably higher due to less efficient operations and highly fragmented networks of small and medium-sized carriers.
%%% Previously separated
Building on our ongoing effort to develop a distributed and dynamic logistics hub network system grounded in the Physical Internet (PI) principles of modular containers and open resource sharing, this study focuses specifically on inter-hub and in-hub operations, with cooperation among multiple regional hub networks. 
Therefore, a shipment may traverse multiple cooperating hub networks. 
With respect to each hub network each shipment enters, it is defined by its expected arrival time at the entry hub and its latest arrival time at the exit hub. 
%%% Previously separated
Based on the defined shipment information, we design a set of multi-hub operation planning protocols for distributed hub operators. 
In their operating networks, operators use our smartly designed protocol separately to plan in-hub shipments' assignments to destination-specific trailers and inter-hub trailers' dispatch schedules.
With carefully designed interconnections between hub networks, the aggregated hub network system is well-positioned to achieve cooperative outcomes and fulfill shipment requests. 
We evaluate the effectiveness of the proposed protocol through a simulation-based experiment under multiple scenarios in an operator's multi-hub network. Overall, this research improves the practicality and robustness of PI-based networks and supports greater cooperation among hub networks toward more efficient and sustainable logistics systems.}
\end{abstract}

\section*{Keywords}
Less-than-Truckload (LTL) Shipping, Physical Internet, Hub Network Operation, Cooperation, Network Planning

\section{Introduction}
Freight flows in modern supply chain systems are characterized by small shipment sizes and high frequency, accelerating reliance on less-than-truckload (LTL) transportation as the primary mode for serving fragmented demand.
Yet, the current LTL industry faces many operational difficulties, especially in balancing efficient consolidation and fast delivery.
These challenges are amplified by the current operational landscape of the LTL industry, where many carriers operate regionally constrained networks, limiting cross-regional consolidation opportunities.

In response, there is growing interest in transportation paradigms to encourage carrier cooperation that shares accessibility of shipments of different carriers at a specific carrier's hubs, allowing freight to traverse interconnected networks. 
Evidence of this shift can be observed in carrier alliances (e.g., Southeastern Freight Lines, A. Duie Pyle, Oak Harbor Freight Lines, and Dayton Freight, in the U.S.) and logistics coordination platforms (e.g, Flexport and DB Schenker). 

These developments align with the principles of the Physical Internet (PI), which envisions standardized, modular, and hyperconnected logistics networks, offering a compelling direction for addressing contemporary challenges in LTL transportation. 
Our ongoing effort is to develop a PI-rooted operating system for LTL shipments stored in PI modular containers built on a hyperconnected logistics network, formed by the cooperation of multiple network operators. 

This study focuses on the hub operation planning protocol for loading shipments and scheduling transportation services (i.e., trucks) to dispatch trailers for each network operator operating a regional hub network. 
Our protocol highlights its practicality, capacity-and-congestion awareness, and ability to interconnect multiple involved operators.
On the planning horizon, the model determines an operator's in-hub shipments' deferral, resorting, and assignment to destination-specific trailers based on their penalties associated with early and late arrivals.
It also plans inter-hub trailers' unloading and dispatch by rigorously considering trailer capacities and hub congestion. 
A multi-scenario simulation-based experiment on an operator's multi-hub network demonstrates the practical feasibility and potential of our protocol in terms of shipment punctuality, trailer fill rates, and scheduled dispatch counts.

The rest of the paper is structured as follows. Section \ref{sec:2} presents the related literature. Section \ref{sec:3} proposes our hub design and protocol. Section \ref{sec:4} analyzes the protocol's effectiveness through a multi-scenario simulation-based experiment. Section \ref{sec:5} summarizes the contributions, limitations, and future work directions.

\section{Literature Review}\label{sec:2}
Early studies on horizontal cooperation in transportation and logistics emphasize its potential to reduce costs, improve environmental performance, and enhance competitiveness, especially for small and medium-sized regional carriers \cite{cbdfs01, b01, aldr01}. However, much of the traditional literature assumes relatively static planning environments, centralized coordination structures, with limited operational integration and restricted scalability \cite{gh01}. To address these limitations, recent research has increasingly adopted the PI into designing operating systems, which envisions open, standardized, and hyperconnected logistics networks \cite{m01, llm01, zlm01, zngma01}. Yet these studies often rely on centralized planning, which is less practical.

Distributed and dynamic planning in PI depends on well-designed inter-hub protocols that govern routing, consolidation, and dispatch \cite{sbphm01}. Shaikh et al. propose a protocol framework for inter-hub transportation services in the PI that enables reliable and scalable hub-to-hub coordination \cite{smhg01}. Shaikh, Muthukrishnan, and Montreuil introduce integrated dynamic freight routing and dispatch protocols that jointly make consolidation and routing decisions \cite{smm01}. 

Nevertheless, past studies on hyperconnected hub network protocols still face practical limitations.
First, they lack rigorous consideration of practical complexities in trailer-loading management. They often assume the hub's loading floor can organize an unconstrained number of loading trailers. In reality, hubs' loading floors only accommodate a limited number of loading trailers (i.e., one per lane) and open a new one only after the previous one is closed. Second, most current network planning research penalizes only late arrivals, while unexpected earliness may also affect overall system quality. Our study addresses these points to fill the gap between the model and real-world practice.

\section{Protocol Design}\label{sec:3}
In this section, we describe the detailed protocol. 
As input, it has knowledge of shipment information, described in subsection~\ref{sec:3.2}, and current hub statuses. 
The protocol's goal is (1) to be capacity-and-congestion-aware for practical and efficient hub operation and (2) to smartly plan consolidation and dispatch of destination-specific trailers through the network.
As output, we obtain the schedules for transportation services to dispatch the trailers.

In subsection~\ref{sec:3.1}, we demonstrate the hub design and the overall protocol. Subsections~\ref{sec:3.2} and~\ref{sec:3.3} describe the detailed protocol that governs shipment loading and inter-hub trailer dispatching, respectively. 

\subsection{Hub Design and Overall Protocol}\label{sec:3.1}
\begin{figure}[htb]
	\centering
	\includegraphics[width=0.8\linewidth]{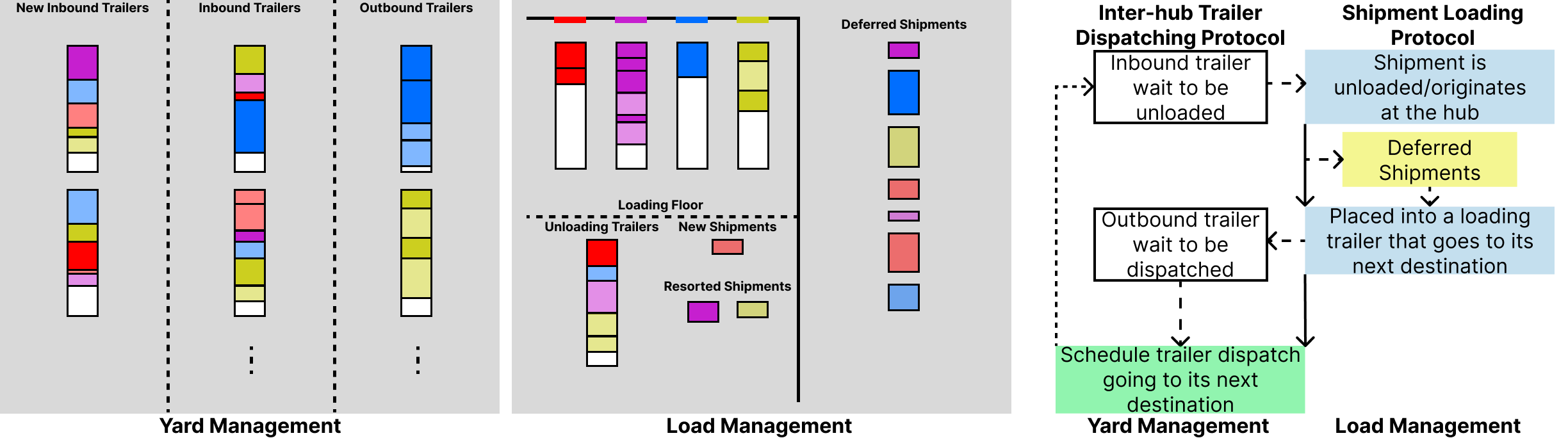}
	\caption{Conceptual Hub Design and Overall Protocol}\label{fig:1}
\end{figure}

Our hub design is inspired by PI hyperconnected hubs \cite{mmb01}. 
We assume hubs across all cooperating network operators' networks share the same design, enabling smooth inter-hub connections, and thus focus on a single hub design and protocol. Also, all shipments are stored in PI modular containers, and another agent plans and supplies sufficient empty trailers in and between hubs.
Figure~\ref{fig:1} shows a conceptual diagram of two interconnected parts in our hub design: load management (\textit{middle}) corresponding to subsection~\ref{sec:3.2}, and yard management (\textit{left}) corresponding to subsection~\ref{sec:3.2}. 

Each colored rectangle represents a shipment, while a tall white rectangle in the background represents a destination specific and one-per-destination trailer. 
The colors (e.g., blue, red, yellow, and purple) indicate shipments' next destination hubs, while the intensities indicate urgencies of leaving the current hub. 
The sizes of colored rectangles indicate the number of modular containers in each shipment. The flowchart on the right demonstrates the overall protocol.

\subsection{Shipment Loading Protocol}\label{sec:3.2}
We denote the network operator's operating hub set by $H$. We denote the estimated travel time from hub $a\in H$ to $b \in H, b\neq a$, by $\tau_{a,b}$ and the estimated handling time at hub $a$ by $h_a$. 
We assume the handling time includes the unloading time and is constant for any number of unloading trailers.
As input, each shipment $s$ has a fixed physical routing of hubs $\mathcal{R}_s$ and arcs $\mathcal{A}_s$ with entry hub $o_s \in H$ and exit hub $d_s \in H$, associated with the estimated entry time $r_s$ and latest exit time $l_s$, respectively. 
Each unsplittable shipment $s$ has $q_s$ modular containers. For shipment $s$ arriving at hub $i$ at planning time point $t_{s, i}$, with remaining routes of nodes $\mathcal{R}^+_s$ and arcs $\mathcal{A}^+_s$, $s$'s estimated total slackness is:
\vspace{0pt}
\begin{equation}
ETS(s, i) = l_s - t_{s, i} - \sum_{a\in\mathcal{R}^+_s}h_a-  \sum_{(a, b) \in\mathcal{A}^+_s} \tau_{a,b}
\label{eq1}
\end{equation}

Each network operator can decide on its unique way to distribute shipment delivery slackness among the remaining hubs. In the study, we assume $ETS(s, i)$ is always equally distributed, thus $s$'s latest departure time from $i$ is:
\vspace{0pt}
\begin{equation}
LDT(s, i) = t_{s, i} + (|R^+_s|)^{-1}\cdot{ETS(s, i)}
\label{eq1}
\end{equation}

The hub $i$ serves a next destination hub set $N_i\subset H$. 
A single destination-specific loading trailer $k_a$ at hub $i$ loads for shipments going to the next hub $a \in N_i$ and has maximum quantity $q^{max}$. 
As the trailers are gradually filled as loading decisions are made, we have $LDT(k_a) = min_{s\in k_a}(LDT(s, i))$ denoting the latest departure time of loading trailer $k_a$ at time point $t$ and $Q(k_a, t) = \sum_{s\in k_a}q_s$ denoting the used space of loading trailer $k_a$ at time point $t$.

\begin{figure}[htb]
	\centering
    \begin{minipage}{0.36\textwidth}
      \centering
      \includegraphics[width=\linewidth]{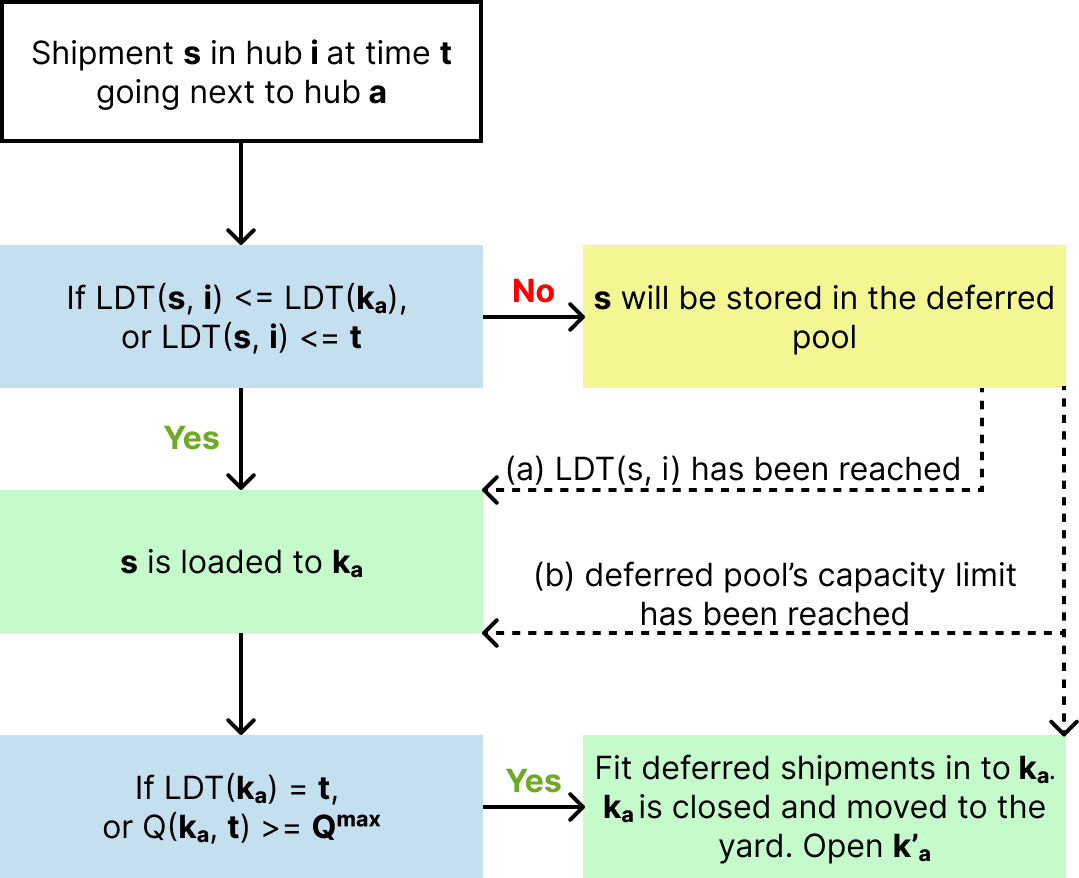}
	\caption{Shipment Loading Protocol}\label{fig:2}
    \end{minipage}\hfill
    \begin{minipage}{0.55\textwidth}
      \centering
      \includegraphics[width=\linewidth]{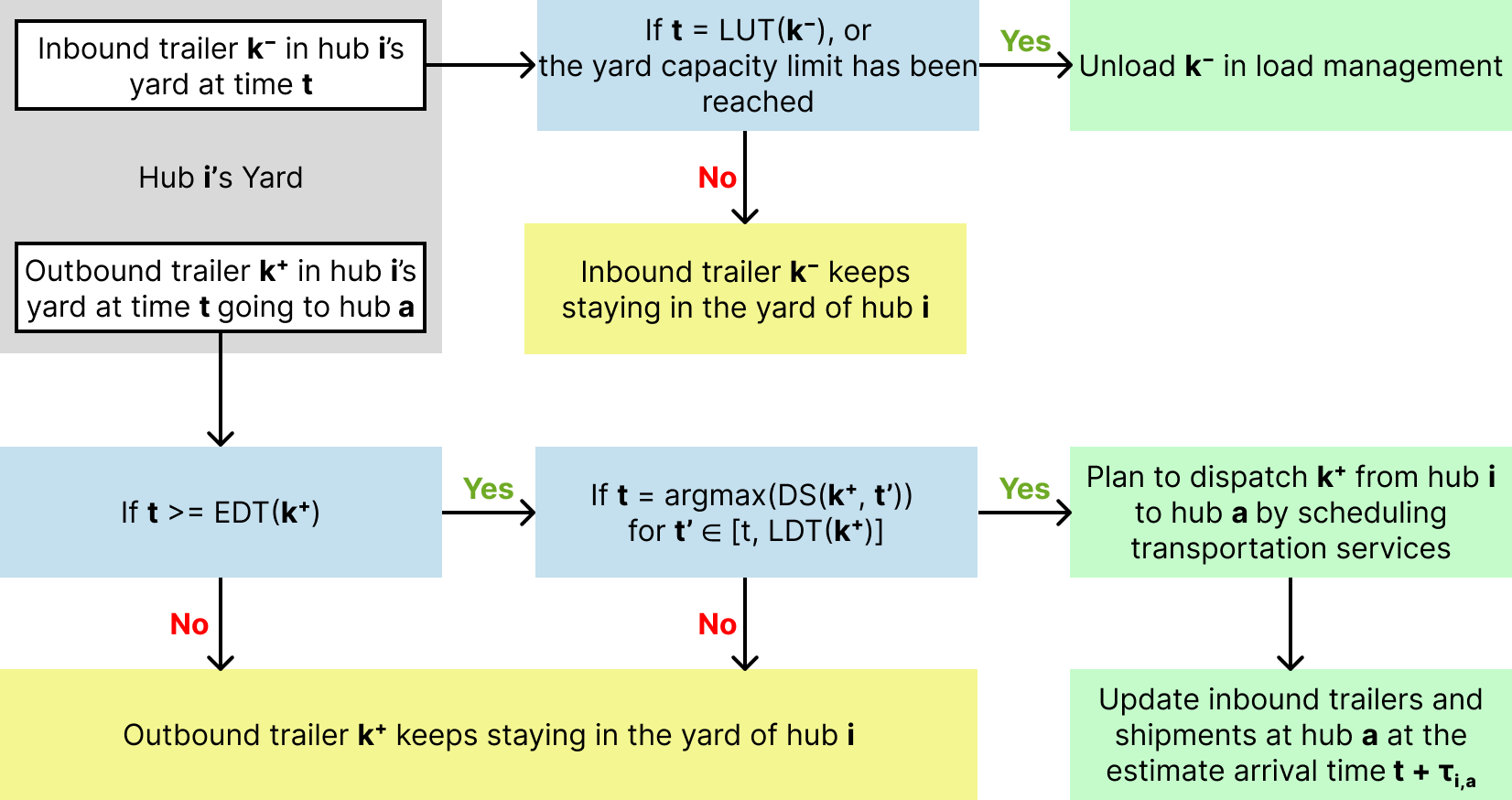}
	\caption{Inter-hub Trailer Dispatching Protocol}\label{fig:3}
    \end{minipage}\hfill
\end{figure}
\vspace{-5pt}

Figure~\ref{fig:2} demonstrates the flow chart for shipment loading protocol in load management. For shipment $s$ unloaded in hub $i$ at time $t$ going next to hub $a$, it will be directly loaded to $k_a$ if $LDT(s, i)$ has been reached or is prior to $LDT(k_a)$. Otherwise, it will be stored in the deferred shipments pool to be resorted later once $LDT(s, i)$ is reached, avoiding delivering non-urgent shipments too early, which incurs penalties. If the deferred pool reaches its capacity limit, shipments in it will be resorted. $k_a$ will be closed and moved to yard if $LDT(k_a, t)$ has been reached or $Q(k_a, t)$ is above a network-operator-preferred threshold $Q^{max}$, and an empty $k_a'$ will be open. Before closing a trailer, we try to fit deferred shipments into it to improve fill rates.

\subsection{Inter-hub Trailer Dispatching Protocol}\label{sec:3.3}
Let $k^-$ and $k^+$ correspondingly denote inbound and outbound trailers. Function $LUT(k^-) = \min_{s\in k^-}(LDT(s, i) - h_i)$ denotes the latest unloading time of an inbound trailer. Function $EDT(k^+)$ denotes the network-operator-preferred earliest departure time of trailer $k^+$, assuming $EDT(k^+) = LDT(k^+) - 6$ in this study.
We define a function $DS(k^+, t')$ that denotes the \textit{dispatch score} of dispatching trailer $k^+$ at future time $t' \in \{t, .., LDT(k^+)\}$. It is calculated as follows:

\vspace{-7pt}
\begin{equation}
DS(k^+, t') = \alpha \cdot A(a, t' + \tau_{ia})
+ \beta \cdot \sum_{s\in k_+} (|k^+|)^{-1} \cdot {U(s, i, t')}
\label{eq1}
\end{equation}
where $\alpha$ and $\beta$ are user-defined coefficients, and $a$ is destination hub of trailer $k^+$. $A(a, t' + \tau_{ia})$ is the estimated destination hub space availability score calculated based on $\sum_{s\in S_{a, t' + \tau_{ia}}} q_s$ versus the hub available space, where $S_{a, t' + \tau_{ia}}$ denotes the set of shipments scheduled to arrive at hub $a$ at $t' + \tau_{ia}$. $U(s, i, t')$ is the utility score of shipment $s$ leaving hub $i$ at $t'$, which is determined by the shipment's early and late penalty. Since the set of shipments scheduled to arrive at the hub at a certain time is dynamically changing, $DS(k^+, t')$ must be iteratively updated every planning time point.

Figure~\ref{fig:3} demonstrates the flowchart for the inter-hub trailer dispatching protocol in yard management. For an inbound trailer $k^-$ awaiting to be unloaded, we iteratively check if its latest unloading time $LUT(k^-)$ has been reached, and unload such a trailer in load management. 
For an outbound trailer $k^+$ waiting to be dispatched in hub $i$'s yard at time $t$ going to hub $a$, we evaluate dispatch scores $DS(k^+, t')$ of time points $t' \in [t, LDT(k^+)]$, if its earliest departure time $EDT(k^+)$ has been reached. If the current time $t = \arg\max_{t'\in[t, LDT(k^+)]}(DS(k^+, t'))$, our protocol plans to dispatch $k^+$ from hub $i$ to hub $a$ at $t$ by scheduling corresponding transportation services and update accordingly.

\section{Results and Discussion}\label{sec:4}
In this study, we create multi-scenario simulation-based experiments to analyze the effectiveness of our protocols within a single hub network operator's multi-hub network. 
We generate 100 distinct cases for each scenario, with 1000 shipment requests per case. 
In each case, we first use the planning protocol to decide trailer dispatches and schedule transportation services based on the scenario setting, and then we evaluate the schedule in a specific simulation experiment generated for the scenario. We consider the following 3 LTL scenarios:

\begin{itemize}
    \setlength{\itemsep}{0pt}
    \setlength{\parskip}{0pt}
    \setlength{\parsep}{0pt}
    \item Scenario 1 assumes hubs with no deferred shipment pool or yard capacity limit, with no uncertainty. \\ (explore the feasible capacity requirements for our protocol's best performance)
    \item Scenario 2 assumes hubs with fixed deferred shipment pool and yard capacity, with no uncertainty. \\ 
    (capacity of deferred shipment pools = 200 modular containers / capacity of yards = 5 trailers)
    \item Scenario 3 assumes hubs with fixed deferred shipment pool and yard capacity, with travel time uncertainty. \\
    (introduce uncertainty $\sim\mathcal{N}(\mu = 0, \sigma=0.4)$ on the arc travel time and plan with the mean travel time)
\end{itemize}

\begin{figure}[htb]
	\centering
	\includegraphics[width=0.8\linewidth]{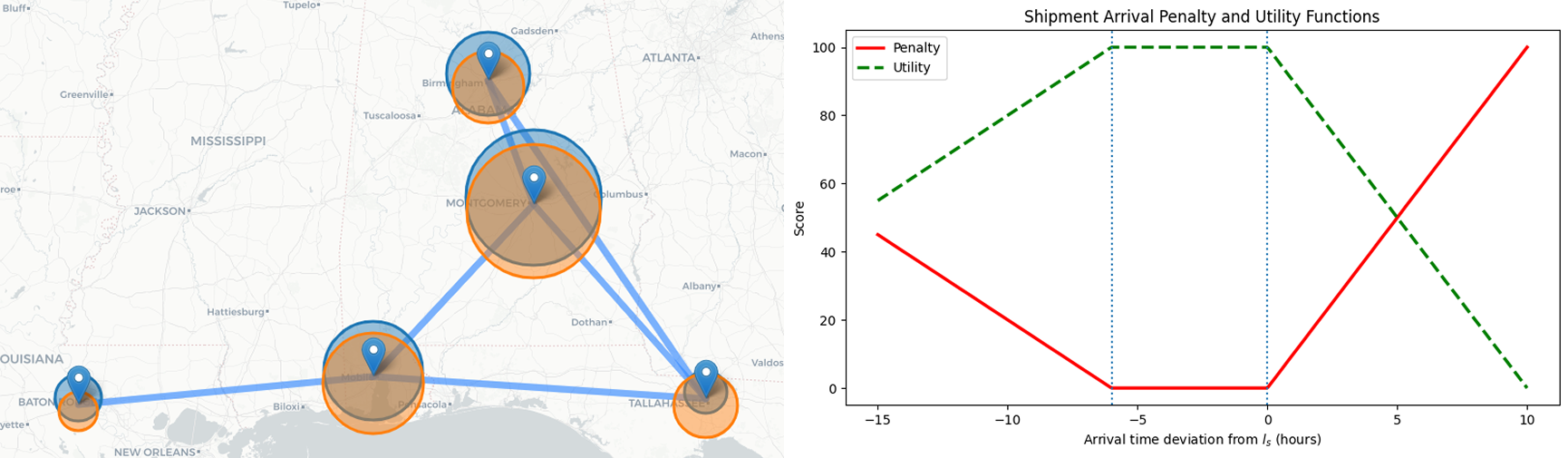}
	\caption{Left: Hub Network and Inbound (Orange) and Outbound (Blue) Shipment Volume; Right: Shipment Arrival Time Penalty (Red) and Utility (Green) Functions}\label{fig:4}
\end{figure}

In evaluating scenario 3, trailers' arrival time might differ from the estimation due to arc travel time uncertainties, incurring deviations from planned loading and dispatch decisions, and potentially leading to the unavailability of transportation services to pick up the trailer.
We report these cases as \textit{failures}.
Figure~\ref{fig:4} visualizes our target regional hub network of five hubs and relay arcs. 
A full-truckload (FTL) scenario is also introduced with no uncertainty. In this scenario, we randomly label five shippers on each shipment, each operating a single trailer for each origin-destination hub pair. FTL shippers ship end-to-end trailers once they are filled or reach the shipments' latest departure time.

For each shipment $s$, we randomly generate its entry hub $o_s$ and exit hub $d_s$ based on each hub's inbound and outbound shipment flows derived from the Freight Analysis Framework \cite{faf}, visualized on the left of Figure~\ref{fig:4}. Each shipment's quantity of modular containers is uniformly randomly selected in $\{1, ..., \frac{q^{max}}{2}\}$, with trailer loading capacity $q^{max} = 50$. The planning time horizon covers 48 hours and is discretized hourly into time points $\{1, ..., 48\}$. We uniformly randomly select entry time $r_s$ in $\{1, ..., 48 - T(o_s, d_s) - RTS(s)\}$, where $T(o_s, d_s)$ denotes the minimum travel time between $o_s$ and $d_s$, and $RTS(s)$ denotes the randomly generated total slackness time for shipment $s$ based on a skewed distribution with mean 12 hours and maximum of 24 hours, and we have exit time $l_s = r_s + T(o_s, d_s) + RTS(s)$.

All shipments share the same penalty score function of a 5-point penalty per hour before $l_s - 6$, a 10-point penalty per hour after $l_s$, and penalty-free for the remainder. The penalty function is used to calculate $U(s, i, t')$ in inter-hub trailer dispatching protocol, as visualized by the right of Figure~\ref{fig:4}. 
Each shipment is assigned the shortest path between its entry and exit hub as the fixed physical routing of hubs $\mathcal{R}_s$ and arcs $\mathcal{A}_s$. We assume full clairvoyance about shipment information, with no difference in shipments between planning and evaluation.

\begin{table}[htb]
\caption{Definitions of Scenarios and Detailed Performances}
\label{tab:1}
\vspace{-0.7cm}
\begin{center}
\begin{tabular}{l|ccc|c}
\hline
 \textbf{Assumptions \& Assessment Metric (a) -- (e)} & \textbf{Scenario 1} & \textbf{Scenario 2} & \textbf{Scenario 3} & \textbf{FTL Scenario}\\ \hline
 \textbf{Assumption (1)} Capacity limits of deferred pool \& yard & $\times$ &  $\circ$  & $\circ$ & - \\
 \textbf{Assumption (2)} Travel time uncertainty & $\times$ &  $\times$  & $\circ$ & $\times$ \\ \hline
\textbf{(a)} Average capacity-enforced deferred shipment resort & - & 617.06 & 534.69 & -\\
\textbf{(b)} Average capacity-enforced inbound trailer unloading & - & 154.32 & 189.55 & -\\
\textbf{(c)} Average capacity-enforced outbound trailer dispatch & - & 102.02 & 75.6 & -\\ \hline
\textbf{(d)} Shipments satisfying the penalty-free window
& 87.02\% & 75.26\% & 72.12\% & 40.28\% \\
\textbf{(e)} Trailers with fill rates $> 90\%$
& 83.13\% & 70.93\% & 66.50\% & 16.45\% \\ 
\textbf{(f)} Average dispatch counts
& 282.49 & 288.85 & 295.94 & 313.43\\\hline
\end{tabular}
\end{center}
\end{table}

\begin{figure}[htb]
	\centering	\includegraphics[width=0.8\linewidth]{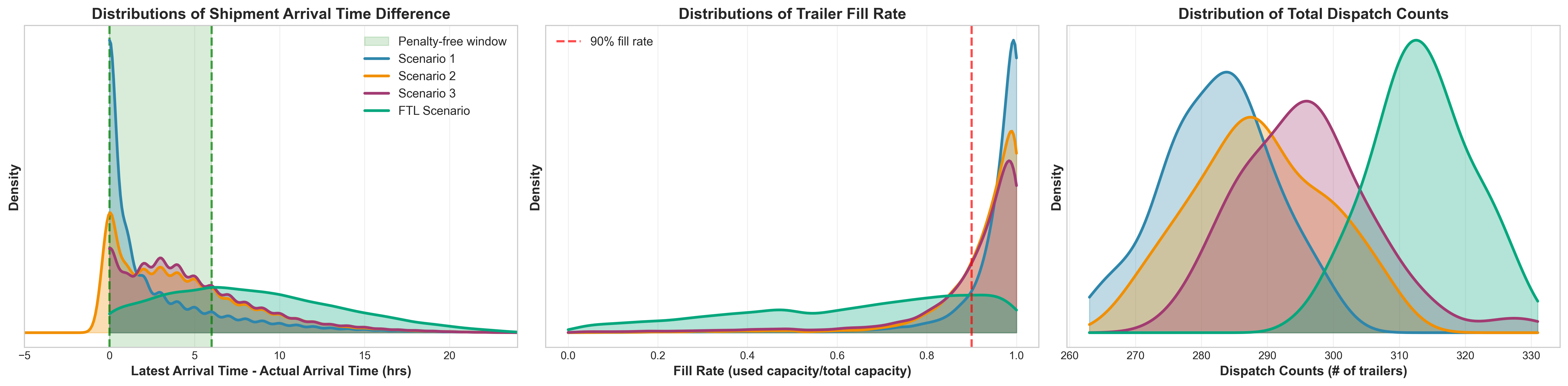}
	\caption{Distributions of  Shipment Arrival Time Difference, Trailer Fill Rate, and Total Dispatch Counts}\label{fig:5}
\end{figure}
Figure~\ref{fig:5} demonstrates our multi-scenario computational results by comparing differences in distributions of shipment arrival time differences (latest arrival time - actual arrival time), trailer fill rates, and total required dispatches. Table~\ref{tab:1} shows the scenario definitions with detailed performances under 6 metrics. 

Scenario 1, with the most relaxed settings, demonstrates the best performances across all metrics (d)--(f). 
The results require a deferred shipment pool capable of storing $1000-1500$ PI modular containers and the yard capable of handling $8-12$ trailers.
In scenario 2, capacity limits might advance time points of inbound trailers' unloading, shipments' resorting, and outbound trailers' dispatch from their desired counterparts: the latest loading time, the latest arrival time, and the dispatch score selected time, respectively. 
The capacity-oriented metrics (a)--(c) in Table~\ref{tab:1} support this interpretation.
Compared to scenario 1, these result in more shipments arriving earlier to the penalty-free window (d), lower trailer fill rates (e), and more average dispatch counts (f). This also indicates needs of hub expansions to accommodate large shipment volumes and a trade-off between infrastructure investment and operational performance. 

Scenario 3 accommodates the strongest practicality among our LTL scenarios at the cost of degraded performances in metrics (d)--(f) compared to other LTL scenarios as shown in Table~\ref{tab:1}.
An average of $75$ (about $\frac{1}{4}$ of scheduled dispatches) failures are observed, incurring extra transportation services for evaluation. 
Despite the trade-off between the performance metric and model realism, the result validates the proposed protocol's viability and motivates refining it to improve resilience in the future.
However, compared to scenario 2, cases of scenario 3 have fewer occurrences of capacity-enforced deferred shipment resort and trailer dispatch.
We attribute this phenomenon to delayed shipments due to travel time uncertainty, which require immediate loading and cannot be deferred.
Trailers containing them also require immediate dispatch, often with low fill rates; they are filled with more deferred shipments before their closing.
Eventually, this cycle leads to lower numbers in metrics (a) and (c) than scenario 2's. 
This observation provides evidence of our protocols' ability to avoid hub congestion under uncertainty, indicating their practical responsiveness and adaptiveness.
Detailed analyses with metrics monitoring hub congestion will be conducted for future studies.

To place our results in a broader context, we also compare our LTL scenarios with the FTL scenario. 
Figure~\ref{fig:5} shows that all LTL scenarios with the proposed protocol successfully outperform the FTL scenario in terms of punctuality and fill rates even under uncertainty. 
Although we count an indirect FTL dispatch (origin and destination hubs not directly connected by relay arcs) as one dispatch, all three LTL scenarios require fewer dispatch counts. This proves our proposed LTL system's success in its fundamental goal of efficiently providing high-quality delivery services and its potential to promote fluid and dynamic shipping behaviors.

\section{Conclusion}\label{sec:5}
This study contributes to the design of distributed operating protocols for hyperconnected LTL operators across multiple regional hub networks. By explicitly modeling practical shipment arrival penalties and hub operations with realistic capacity and congestion constraints, our work bridges an important gap between PI-inspired conceptual frameworks and the operational realities faced by LTL carriers. 
We propose a set of distributed protocols that govern shipment loading, trailer consolidation, and inter-hub dispatch. 
Through multi-scenario, simulation-based experiments, we show that although incorporating realistic operational constraints inevitably leads to some performance degradation, the proposed approach maintains a competitive level of shipment punctuality and trailer utilization, demonstrating its practical viability.
Future research will extend this work to incorporate rolling-horizon planning, monitor hub congestion, and conduct multi-operator experiments to provide stronger support to cooperation among network operators.

\section*{Acknowledgements}
We thank Praveen Muthukrishnan from Georgia Institute of Technology's Physical Internet Center for sharing data and knowledge on hyperconnected networks.

\end{document}